# Approximation simultanée d'un nombre et de son carré

**Damien Roy**

**Département de mathématiques, Université d'Ottawa, 585 King Edward, Ottawa, Ontario K1N 6N5, Canada**
**Courriel : droy@uottawa.ca**

**Résumé.** En 1969, H. Davenport et W. Schmidt ont établi une mesure d'approximation simultanée pour un nombre réel $\xi$ et son carré par des nombres rationnels de même dénominateur. Cette mesure suppose seulement que le nombre $\xi$ n'est ni rationnel ni quadratique sur $\mathbf{Q}$. Nous montrons ici, par un exemple, que cette mesure est optimale. Nous indiquons aussi plusieurs propriétés des nombres pour lesquels cette mesure est optimale, notamment en regard à l'approximation par les entiers algébriques de degré au plus trois. © 2001 Académie des sciences/Éditions scientifiques et médicales Elsevier SAS

*Simultaneous approximation to a real number and its square*

**Abstract.** *In 1969, H. Davenport and W. Schmidt established a measure of simultaneous approximation for a real number $\xi$ and its square by rational numbers with the same denominator, assuming only that $\xi$ is not rational nor quadratic over $\mathbf{Q}$. Here, we show by an example, that this measure is optimal. We also indicate several properties of the numbers for which this measure is optimal, in particular with respect to approximation by algebraic integers of degree at most three.* © 2001 Académie des sciences/Éditions scientifiques et médicales Elsevier SAS

## *Abridged English version*

### 1. Introduction.

Let $\xi$ be a real number and let $\varphi\colon [1,\infty) \to (0,\infty)$ be a non-increasing function. We look for integer triples $\mathbf{x} = (x_0, x_1, x_2) \in \mathbf{Z}^3$ which, for a given real number $X \geqslant 1$, satisfy the conditions

$$0 < x_0 \leqslant X, \quad |x_0\xi - x_1| \leqslant \varphi(X) \quad \text{and} \quad |x_0\xi^2 - x_2| \leqslant \varphi(X). \tag{1.1}$$

When such a triple exists, the quotients $x_1/x_0$ and $x_2/x_0$ provide approximations to $\xi$ and $\xi^2$ by rational numbers with the same denominator $x_0$. We recall the following facts.

---

                                         1

**Damien Roy**

(1) Dirichlet box principle shows that (1.1) admits an integral solution for each $X \geqslant 1$ by taking $\varphi(X) = [\sqrt{X}]^{-1}$ where $[\sqrt{X}]$ denotes the integral part of $\sqrt{X}$ (see for example Theorem 1A in Chapter II of [7]).

(2) If $\xi$ is algebraic over $\mathbf{Q}$ of degree $\leqslant 2$, the same is true with the function $\varphi(X) = cX^{-1}$ for an appropriate constant $c > 0$, depending on $\xi$.

(3) Conversely, Davenport and Schmidt proved that, for any real number $\xi$ with $[\mathbf{Q}(\xi)\colon \mathbf{Q}] > 2$, there exists a constant $c = c(\xi) > 0$ such that, for arbitrarily large values of $X$, the inequalities (1.1) admit no integral solution with $\varphi(X) = cX^{-1/\gamma}$ where $\gamma = (1+\sqrt{5})/2$ denotes the golden ratio (Theorem 1a of [4]). Numerically, we have $1/\gamma \simeq 0.618$.

(4) For $\varphi(X) = X^{-\lambda}$ with $\lambda > 1/2$, the set of real numbers $\xi$ for which (1.1) admit integral solutions for arbitrarily large values of $X$ has Lebesgue measure zero. This follows as in §4 of [2] from a metrical result of Sprindžuk via Lemma 1 of [2] and a duality argument (see page 395 of [4]). By the subspace theorem of Schmidt, the same set contains no algebraic number of degree at least three (see for example Theorem 1B in Chapter VI of [7]).

In [6], we show that the result of Davenport and Schmidt mentioned in (3) is optimal up to the value of the constant:

THEOREM 1.1. – *There exists a real number $\xi$ which is neither rational nor quadratic irrational and which has the property that, for a suitable constant $c > 0$, the inequalities*

$$0 < x_0 \leqslant X, \quad |x_0\xi - x_1| \leqslant cX^{-1/\gamma} \quad et \quad |x_0\xi^2 - x_2| \leqslant cX^{-1/\gamma} \tag{1.2}$$

*have a non-zero solution $(x_0, x_1, x_2) \in \mathbf{Z}^3$ for any real number $X \geqslant 1$. Any such number is transcendental over $\mathbf{Q}$ and the set of these real numbers is enumerable.*

The fact that, for such a number, (1.2) admits an integral solution for each $X \geqslant 1$ is a notable improvement over the estimates provided by the box principle and stated in item (1) above, since $1/\gamma > 1/2$. Such a phenomenon reminds of Cassels' counter-example in the context of algebraic independence (see Theorem XIV in Chapter V of [3]) but it is surprising to see it happen in a problem of approximation to numbers ($\xi$ and $\xi^2$) belonging to a field of transcendence degree one over $\mathbf{Q}$.

Examples of such numbers are the real numbers

$$\xi_{a,b} = [0, a, b, a, a, b, a, \ldots] = 1/(a + 1/(b + \ldots))$$

whose continued fraction expansion is given by the Fibonacci word on two distinct positive integers $a$ and $b$, that is the infinite word $abaabab\ldots$ starting with $a$ which is the fixed point of the substitution sending $a$ to $ab$ and $b$ to $a$. They constitute a special case of the more general Sturmian continued fractions studied by Allouche, Davison, Queffélec and Zamboni in [1]. These authors proved that the latter are transcendental over $\mathbf{Q}$ by showing that they admit very good approximations by quadratic real numbers, which accounts for the fact that any Sturmian sequence begins in arbitrarily long squares (see Proposition 2 of [1]). Here, we use a different combinatorial property of the Fibonacci word, namely the fact that, if $f(i)$ denotes the $i$-th Fibonacci number defined recursively by $f(0) = f(1) = 1$ and $f(i) = f(i-1) + f(i-2)$ for $i \geqslant 2$, then, for $i \geqslant 3$, the first $f(i) - 2$ letters of the Fibonacci word form a palindrome. From this, we find that the numbers $\xi_{a,b}$ satisfy the following property, stronger than that required in Theorem 1.1 (see the French version for a proof).

THEOREM 1.2. – *Let $a$, $b$ be distinct positive integers. There exist positive constants $c_1, c_2, c_3$ and a sequence $(\mathbf{x}_i)_{i \geqslant 1}$ of integer triples with non-negative entries such that, putting $\mathbf{x}_i = (x_{i,0}, x_{i,1}, x_{i,2})$ and $X_i = x_{i,0}$, we have, for $i \geqslant 2$,*

$$c_1 X_{i-1}^\gamma < X_i < c_2 X_{i-1}^\gamma \quad \text{and} \quad \max_{j=1,2} \left| x_{i,0} \xi_{a,b}^j - x_{i,j} \right| \leqslant c_3 X_i^{-1}.$$

Davenport and Schmidt deduce from their result stated in item (3) above that, for any real number $\xi$ with $[\mathbf{Q}(\xi) : \mathbf{Q}] \geqslant 2$, there exists another constant $c > 0$ and infinitely many algebraic integers $\alpha$ of degree $\leqslant 3$ which satisfy $|\xi - \alpha| \leqslant c H(\alpha)^{-\gamma^2}$, where $H(\alpha)$ stands for the usual height of $\alpha$, that is the largest absolute value of the coefficients of its irreducible polynomial over $\mathbf{Z}$ (Theorem 1 of [4]).

Conversely, one may hope that, for any real number $\theta > \gamma^2$ and any $\xi$ as in Theorem 1.1, there exist only finitely many algebraic integers $\alpha$ of degree $\leqslant 3$ which satisfy $|\xi - \alpha| \leqslant X^{-\theta}$. To show this for the Fibonacci continued fractions $\xi_{a,b}$, it would suffice to prove, in the notation of Theorem 1.2, that, for any $\delta > 0$, there exist only finitely many indices $i$ for which the distance from $x_{i,0} \xi_{a,b}^3$ to the nearest integer is $\leqslant X^{-\delta}$ (see Proposition 9.1 of [6]). Numerical experiments seem to be in agreement with this. The result below presents a partial step in this direction as well as other approximation properties of these numbers (see §§7–9 of [6]):

THEOREM 1.3. – *Let $\xi$ be as in Theorem 1.1. Their exist positive constants $c_1, \ldots, c_5$ with the following properties:*

*(1) for any rational number $\alpha \in \mathbf{Q}$, we have $|\xi - \alpha| \geqslant c_1 H(\alpha)^{-2}(1 + \log H(\alpha))^{-c_2}$;*

*(2) for any algebraic number $\alpha \in \mathbf{C}$ of degree at most two over $\mathbf{Q}$, we have $|\xi - \alpha| \geqslant c_3 H(\alpha)^{-2\gamma^2}$;*

*(3) there exist infinitely many quadratic real numbers $\alpha$ with $|\xi - \alpha| \leqslant c_4 H(\alpha)^{-2\gamma^2}$;*

*(4) for any algebraic integer $\alpha \in \mathbf{C}$ of degree at most three over $\mathbf{Q}$, we have $|\xi - \alpha| \geqslant c_5 H(\alpha)^{-(3/2)\gamma^2}$.*

---

## 1. Introduction.

Soit $\xi$ un nombre réel et soit $\varphi \colon [1, \infty) \to (0, \infty)$ une fonction décroissante. On s'intéresse aux triplets d'entiers $\mathbf{x} = (x_0, x_1, x_2) \in \mathbf{Z}^3$ qui, pour un réel $X \geqslant 1$ donné, satisfont les conditions

$$0 < x_0 \leqslant X, \quad |x_0 \xi - x_1| \leqslant \varphi(X) \quad \text{et} \quad |x_0 \xi^2 - x_2| \leqslant \varphi(X). \tag{1.1}$$

Lorsqu'un tel triplet existe, les quotients $x_1/x_0$ et $x_2/x_0$ constituent des approximations de $\xi$ et $\xi^2$ par des nombres rationnels de même dénominateur $x_0$. On rappelle les faits suivants.

(1) Le principe des tiroirs de Dirichlet montre que (1.1) admet une solution entière pour tout $X \geqslant 1$ avec la fonction $\varphi(X) = [\sqrt{X}]^{-1}$ où $[\sqrt{X}]$ désigne la partie entière de $\sqrt{X}$ (voir par exemple le théorème 1A du chapitre II de [7]).

(2) Si $\xi$ est algébrique sur $\mathbf{Q}$ de degré $\leqslant 2$, c'est encore vrai avec la fonction $\varphi(X) = cX^{-1}$ pour une constante $c > 0$ appropriée, fonction de $\xi$.

**Damien Roy**

(3) En contrepartie, Davenport et Schmidt ont montré que, pour tout nombre réel $\xi$ avec $[\mathbf{Q}(\xi)\colon \mathbf{Q}] > 2$, il existe une constante $c = c(\xi) > 0$ telle que, pour des valeurs de $X$ arbitrairement grandes, les inégalités (1.1) n'admettent pas de solution entière avec $\varphi(X) = cX^{-1/\gamma}$ où $\gamma = (1+\sqrt{5})/2$ désigne le nombre d'or (Théorème 1a de [4]). Numériquement, on a $1/\gamma \simeq 0.618$.

(4) Pour $\varphi(X) = X^{-\lambda}$ avec $\lambda > 1/2$ fixé, l'ensemble des nombres réels $\xi$ pour lesquels les inégalités (1.1) admettent des solutions entières pour des valeurs arbitrairement grandes de $X$ est un ensemble de mesure de Lebesgue nulle. Cela découle d'un résultat métrique de Sprindžuk via le lemme 1 de [2] et un argument de dualité (voir page 395 de [4]). En vertu du théorème du sous-espace de Schmidt, le même ensemble ne contient aucun nombre réel algébrique sur $\mathbf{Q}$ de degré au moins trois (voir par exemple le théorème 1B du chapitre VI de [7]).

Dans [6], on montre que le résultat de Davenport et Schmidt indiqué en (3) est optimal à la valeur de la constante près :

THÉORÈME 1.1. – *Il existe un nombre réel $\xi$ qui est transcendant sur $\mathbf{Q}$ et, pour ce nombre, une constante $c = c(\xi) > 0$ telle que, pour tout $X \geqslant 1$, les inégalités*

$$0 < x_0 \leqslant X, \quad |x_0\xi - x_1| \leqslant cX^{-1/\gamma} \quad et \quad |x_0\xi^2 - x_2| \leqslant cX^{-1/\gamma} \tag{1.2}$$

*admettent une solution en entiers $x_0, x_1, x_2$. L'ensemble de ces nombres réels $\xi$ est dénombrable.*

Le fait que, pour un tel nombre, (1.2) admette une solution entière pour tout $X \geqslant 1$ améliore nettement les bornes fournies par le principe des tiroirs et décrites ci-dessus en (1), puisque $1/\gamma > 1/2$. Un tel phénomène rappelle la situation du contre-exemple de Cassels en indépendance algébrique (voir le théorème XIV du chapitre V de [3]) mais il surprenant dans le voir surgir dans des questions d'approximation de nombres ($\xi$ et $\xi^2$) appartenant à un corps de degré de transcendance un sur $\mathbf{Q}$.

Le lecteur pourra consulter [6] pour une étude systématique de ces nombres. Nous nous contenterons ici d'en donner un exemple, les fractions continues de Fibonacci. Nous mentionnerons aussi d'autres propriétés de ces nombres reliées à l'approximation par les nombres algébriques et tirées de [6].

## 2. Fractions continues de Fibonacci

Soit $E = \{a, b\}$ un alphabet de deux lettres et soit $E^*$ le monoïde des mots sur $E$ avec pour produit la concaténation des mots. La *suite de Fibonacci* dans $E^*$ est la suite de mots $(w_i)_{i \geqslant 0}$ definie récursivement en posant

$$w_0 = b, \quad w_1 = a \quad \text{et} \quad w_i = w_{i-1}w_{i-2}, \quad (i \geqslant 2)$$

(voir l'exemple 1.3.6 de [5]). Puisque, pour tout $i \geqslant 1$, le mot $w_i$ est un préfixe de $w_{i+1}$, cette suite converge vers un mot infini $w = abaabaab\ldots$ appelé *mot de Fibonacci* sur $\{a, b\}$.

Soient $a$ and $b$ des entiers positifs distincts et soit

$$\xi = \xi_{a,b} = [0, w] = [0, a, b, a, a, b, a, \ldots] = 1/(a + 1/(b + \ldots))$$

le nombre réel dont la suite des quotients partiels est $0$ suivi par le mot de Fibonacci sur $\{a, b\}$. C'est un exemple de fraction continue sturmienne au sens d'Allouche, Davison, Queffélec et Zamboni [1]. Ces auteurs ont montré que de tels nombres sont transcendants sur $\mathbf{Q}$ en vertu de leur propriété d'approximation par les nombres quadratiques réels qui découle du fait que toute suite sturmienne, sans être ultimement périodique, débute par des carrés de longueur arbitrairement grande (Proposition 2 de [1]). Nous utiliserons ici une autre propriété du mot de Fibonacci :

LEMME 2.1. – *Pour tout entier $i \geqslant 1$, le mot $w_{i+2}$ privé de ses deux dernières lettres est un palindrome $m_i$. De plus, si on pose $s_i = ab$ pour $i$ pair et $s_i = ba$ pour $i$ impair, on a les formules de récurrence*

$$m_1 = a, \quad m_2 = aba \quad et \quad m_i = m_{i-1}s_{i-1}m_{i-2}, \quad (i \geqslant 3). \tag{2.1}$$

*Démonstration.* Les formules (2.1) se vérifient sans peine en observant que $w_i$ se termine par $s_i$ pour tout $i \geqslant 2$. On vérifie aussi que $m_i$ est un palindrome pour $i = 1, 2, 3$. Enfin, si on suppose $i \geqslant 4$ et que $m_j$ est un palindrome pour $j = 1, \ldots, i-1$, on trouve

$$m_i = m_{i-1}s_{i-1}m_{i-2} = m_{i-2}s_{i-2}m_{i-3}s_{i-1}m_{i-2}.$$

Donc, $m_i$ est aussi un palindrome. □

THÉORÈME 2.2. – *Soient $a$ et $b$ des entiers positifs distincts et soit $\xi = \xi_{a,b}$. Il existe des constantes positives $c_1, c_2, c_3$ et une suite $(\mathbf{x}_i)_{i \geqslant 1}$ de triplets d'entiers positifs ou nuls telles qu'en posant $\mathbf{x}_i = (x_{i,0}, x_{i,1}, x_{i,2})$ et $X_i = x_{i,0}$, on ait, pour $i \geqslant 2$,*

$$c_1 X_{i-1}^\gamma < X_i < c_2 X_{i-1}^\gamma \quad et \quad \max_{j=1,2} \left| x_{i,0} \xi^j - x_{i,j} \right| \leqslant c_3 X_i^{-1}. \tag{2.2}$$

On en déduit sans peine que, pour tout $X \geqslant X_2$, il existe un triplet d'entiers $\mathbf{x} = (x_0, x_1, x_2)$ qui vérifie les inégalités (1.2) avec $c = c_3 c_2^{1/\gamma}$. Il suffit en effet de choisir $\mathbf{x} = \mathbf{x}_i$ où $i$ est un entier tel que $X_i \leqslant X \leqslant X_{i+1}$.

*Démonstration.* On note d'abord que, si on écrit $\xi = [0, a_1, a_2, \ldots]$, alors les convergents $[0, a_1, \ldots, a_j]$ de $\xi$ s'écrivent sous forme réduite $p_j/q_j$ avec des entiers $p_j \geqslant 0$, $q_j > 0$ donnés par les relations matricielles

$$\begin{pmatrix} q_j & q_{j-1} \\ p_j & p_{j-1} \end{pmatrix} = \begin{pmatrix} a_1 & 1 \\ 1 & 0 \end{pmatrix} \begin{pmatrix} a_2 & 1 \\ 1 & 0 \end{pmatrix} \cdots \begin{pmatrix} a_j & 1 \\ 1 & 0 \end{pmatrix}$$

et de plus on a

$$|q_j \xi - p_j| \leqslant q_j^{-1}$$

(voir le chapitre I de [7]). Pour appliquer ce résultat, on pose

$$A = \begin{pmatrix} a & 1 \\ 1 & 0 \end{pmatrix} \quad \text{et} \quad B = \begin{pmatrix} b & 1 \\ 1 & 0 \end{pmatrix},$$

et on note $\Phi \colon E^* \to \mathrm{GL}_2(\mathbf{Z})$ l'homomorphisme de monoïdes qui applique $a$ sur $A$ et $b$ sur $B$. Comme $A$ et $B$ sont des matrices symétriques et que les mots $m_i$ sont des palindromes, les matrices $\Phi(m_i)$ sont symétriques. On peut donc écrire

$$M_i = \Phi(m_i) = \begin{pmatrix} x_{i,0} & x_{i,1} \\ x_{i,1} & x_{i,2} \end{pmatrix}$$

pour un triplet d'entiers $\mathbf{x}_i = (x_{i,0}, x_{i,1}, x_{i,2})$. On en déduit que $x_{i,0}$ et $x_{i,1}$ sont positifs, que $x_{i,1}/x_{i,0} = [0, m_i]$ est un convergent de $\xi$ et que $x_{i,2}/x_{i,1}$ est le convergent précédent. Donc, on a

$$|x_{i,0} \xi - x_{i,1}| \leqslant x_{i,0}^{-1} \quad \text{et} \quad |x_{i,1} \xi - x_{i,2}| \leqslant x_{i,1}^{-1}$$

et par suite, en posant $X_i = x_{i,0}$, la seconde des conditions (2.2) est satisfaite pour une constante $c_3$.

Pour $i \geqslant 3$, la relation (2.1) livre $M_i = M_{i-1} S_{i-1} M_{i-2}$ avec $S_{i-1} = \Phi(s_{i-1})$, donc

$$x_{i,0} = \begin{pmatrix} x_{i-1,0} & x_{i-1,1} \end{pmatrix} S_{i-1} \begin{pmatrix} x_{i-2,0} \\ x_{i-2,1} \end{pmatrix},$$

et par suite

$$\lim_{i \to \infty} \frac{X_i}{X_{i-2} X_{i-1}} = \begin{pmatrix} 1 & \xi \end{pmatrix} AB \begin{pmatrix} 1 \\ \xi \end{pmatrix} = \xi^2 + (a+b)\xi + (ab+1). \tag{2.3}$$

**Damien Roy**

Posons $q_i = X_i X_{i-1}^{-\gamma}$ pour tout entier $i \geqslant 2$. Comme $\gamma = 1 + 1/\gamma$, on trouve

$$q_i = \frac{X_i}{X_{i-2}X_{i-1}} q_{i-1}^{-1/\gamma}.$$

En vertu de (2.3), il existe donc des constantes $c_5 > c_4 > 0$ telles que

$$c_4 q_{i-1}^{-1/\gamma} < q_i < c_5 q_{i-1}^{-1/\gamma}$$

pour tout entier $i \geqslant 3$. Par récurrence sur $i$, on en déduit $c_1 < q_i < c_2$ pour tout $i \geqslant 2$, en choisissant $c_1 = \min\{q_2, c_4^\gamma/c_5\}$ et $c_2 = \max\{q_2, c_5^\gamma/c_4\}$ et la première des conditions (2.2) est elle-aussi satisfaite. □

## 3. Autres propriétés

Davenport et Schmidt déduisent de leur résultat rappelé à l'item (3) de l'introduction que, pour tout nombre réel $\xi$ avec $[\mathbf{Q}(\xi) : \mathbf{Q}] \geqslant 2$, il existe encore une constante $c > 0$ et une infinité d'entiers algébriques $\alpha$ de degré $\leqslant 3$ qui vérifient $|\xi - \alpha| \leqslant cH(\alpha)^{-\gamma^2}$, où $H(\alpha)$ désigne la hauteur naïve de $\alpha$, c'est-à-dire le maximum des valeurs absolues de son polynôme irréductible sur $\mathbf{Z}$ (Théorème 1 de [4]).

Il est possible que ce second résultat soit lui-aussi essentiellement optimal contrairement à ce qu'il est naturel de conjecturer (voir page 259 de [7]). Plus précisément, on peut espérer que, pour tout réel $\theta > \gamma^2$ et tout nombre $\xi$ comme au théorème 1.1, il n'existe qu'un nombre fini d'entiers algébriques $\alpha$ de degré $\leqslant 3$ qui satisfont $|\xi - \alpha| \leqslant X^{-\theta}$. Pour l'établir pour les fractions continues de Fibonacci $\xi = \xi_{a,b}$, il suffirait de montrer, dans les notations du théorème 2.2, que, pour tout $\delta > 0$, il n'existe qu'un nombre fini d'indices $i$ pour lesquels la distance de $x_{i,0}\xi^3$ au plus proche entier soit $\leqslant X^{-\delta}$ (voir la proposition 9.1 de [6]). Les expériences numériques que nous avons faites militent d'ailleurs en ce sens. Nous savons toutefois montrer un résultat partiel dans cette direction ainsi que d'autres propriétés d'approximation de ces nombres (voir §§7–9 de [6]) :

THÉORÈME 3.1. – *Soit $\xi$ comme au théorème 1.1. Il existe des constantes positives $c_1, \ldots, c_5$ avec les propriétés suivantes :*

*(1) pour tout nombre rationnel $\alpha \in \mathbf{Q}$, on a $|\xi - \alpha| \geqslant c_1 H(\alpha)^{-2}(1 + \log H(\alpha))^{-c_2}$ ;*

*(2) pour tout nombre algébrique $\alpha \in \mathbf{C}$ de degré au plus deux sur $\mathbf{Q}$, on a $|\xi - \alpha| \geqslant c_3 H(\alpha)^{-2\gamma^2}$ ;*

*(3) il existe une infinité de nombres réels quadratiques $\alpha$ avec $|\xi - \alpha| \leqslant c_4 H(\alpha)^{-2\gamma^2}$ ;*

*(4) pour tout entier algébrique $\alpha \in \mathbf{C}$ de degré au plus trois sur $\mathbf{Q}$, on a $|\xi - \alpha| \geqslant c_5 H(\alpha)^{-(3/2)\gamma^2}$.*